\documentclass{article}
\begin{document}

\title{On the Coefficients of a Fibonacci Power Series}
\author{Federico Ardila M.}
\date{October 15, 2001}
\maketitle

\begin{abstract}
We give an explicit description of the coefficients of the formal
power series $(1-x) (1-x^2) (1-x^3) (1-x^5) (1-x^8) (1-x^{13})
\cdots$. In particular, we show that all the coefficients are
equal to $-1, 0$ or $1$.
\end{abstract}

\medskip

The Fibonacci numbers are defined by the recurrence relation
$F_{n+2} = F_{n+1} + F_n$ for $n \geq 0$, and the initial
conditions $F_0=0, F_1=1$. Consider the infinite product
\begin{displaymath}
\begin{array}{rcl}
A(x) & = & \prod_{k \geq 2} (1-x^{F_k}) =
(1-x)(1-x^2)(1-x^3)(1-x^5)(1-x^8)\cdots \\
\\
& = & 1- x - x^2 + x^4 + x^7 - x^8 + x^{11} - x^{12} - x^{13} +
x^{14} + x^{18} + \cdots
\end{array}
\end{displaymath}
regarded as a formal power series. In \cite{4}, N. Robbins proved
that the coefficients of $A(x)$ are all equal to $-1, 0$ or $1$.
We shall give a short proof of this fact, and a very simple
recursive description of the coefficients of $A(x)$.

Following the notation of \cite{4}, let $a(m)$ be the coefficient
of $x^m$ in $A(x)$. It is clear that $a(m)=r_E(m)-r_O(m)$, where
$r_E(m)$ is equal to the number of partitions of $m$ into an even
number of distinct positive Fibonacci numbers, and $r_O(m)$ is
equal to the number of partitions of $m$ into an odd number of
distinct positive Fibonacci numbers. We call these partitions
``even" and ``odd" respectively.

\medskip

{\bf Proposition 1.} Let $n \geq 5$ be an integer. Consider the
coefficients $a(m)$ for $m$ in the interval $[F_n,F_{n+1})$.
Split this interval into the three subintervals $[F_n,
F_n+F_{n-3}-2], [F_n+F_{n-3}-1, F_n+F_{n-2}-1]$ and
$[F_n+F_{n-2}, F_{n+1}-1]$.
\begin{enumerate}
\item
The numbers $a(F_n), a(F_n+1), \ldots, a(F_n+F_{n-3}-2)$ are
equal to the numbers $(-1)^{n-1}a(F_{n-3}-2),
(-1)^{n-1}a(F_{n-3}-3), \ldots, (-1)^{n-1}a(0)$ in that order.
\item
The numbers $a(F_n+F_{n-3}-1), a(F_n+F_{n-3}), \ldots,
a(F_n+F_{n-2}-1)$ are equal to $0$.
\item
The numbers $a(F_n+F_{n-2}), a(F_n+F_{n-2}+1), \ldots,
a(F_{n+1}-1)$ are equal to the numbers $a(0), a(1), \ldots,
a(F_{n-3}-1)$ in that order.
\end{enumerate}

\medskip

This description gives a very fast method for computing the
coefficients $a(m)$ recursively. Once we have computed them for
$0 \leq m < F_n$ we can immediately compute them for $F_n \leq m <
F_{n+1}$ using Proposition 1.

Also, since the coefficient of $x^m$ in $A(x)$ is equal to $-1, 0$
or $1$ for all non-negative integers $m < F_5$, it follows
inductively that the coefficients in each interval $[F_n,
F_{n+1})$ are also all equal to $-1, 0$ or $1$. This proves
Robbins's result.

\medskip

{\bf Proof of Proposition 1.} It will be convenient to prove
Proposition 1.2 first. Let $F_n+F_{n-3}-1 \leq m \leq
F_n+F_{n-2}-1$, and consider the partitions of $m$ into distinct
positive Fibonacci numbers. It is clear that the largest part in
such a partition cannot be $F_{n+1}$ or larger. It cannot be
$F_{n-2}$ or smaller either, because $F_{n-2}+F_{n-3} + \cdots +
F_2 = F_n - 2 < m$. Therefore, it must be $F_n$ or $F_{n-1}$.

If the largest part is $F_n$, then the second largest part cannot
be $F_{n-1}$ or $F_{n-2}$. If, on the other hand, it is
$F_{n-1}$, then the second largest part must be $F_{n-2}$,
because $F_{n-1}+F_{n-3}+F_{n-4}+\cdots+F_2 = 2F_{n-1}-2 = F_n +
F_{n-3}-2 < m$.

This means that we can split the set of partitions into pairs.
Each pair consists of two partitions of the form $F_n + F_a + F_b
+ \cdots$ and $F_{n-1} + F_{n-2} + F_a + F_b + \cdots$, where $n-3
\geq a > b > \ldots$. In each pair, one of the partitions is even
and the other is odd. Therefore $r_E(m)=r_O(m)$ and $a(m)=0$ as
claimed.

\medskip

Now we use a similar analysis to prove Proposition 1.3. Let $F_n
+ F_{n-2} \leq m \leq F_{n+1}-1$. As before, the largest part of
a partition of $m$ must be $F_n$ or $F_{n-1}$. If it is $F_n$,
the second largest part cannot be $F_{n-1}$. If, on the other
hand, it is $F_{n-1}$, then the second largest part must be
$F_{n-2}$.

Again, we can split {\it a subset of the set of partitions} into
pairs. Each pair consists of two partitions of the form $F_n +
F_a + F_b + \cdots$ and $F_{n-1} + F_{n-2} + F_a + F_b + \cdots$,
where $n-3 \geq a > b > \ldots$. In each pair there is an even
and an odd partition.

The remaining partitions are of the form $F_n + F_{n-2} + F_a +
F_b + \cdots$, where $n-3 \geq a > b > \ldots$. To each one of
these partitions we can assign a partition of $m' = m - F_n -
F_{n-2}$, by just removing the parts $F_n$ and $F_{n-2}$. This is
in fact a bijection. Since $m' < F_{n-2}$, any partition of $m'$
has largest part less than or equal to $F_{n-3}$; therefore it
can be obtained in that way from a partition of $m$.

It is clear that, under this bijection, odd partitions of $m$ go
to odd partitions of $m'$ and even partitions of $m$ go to even
partitions of $m'$. It follows that $a(m) = a(m - F_n -
F_{n-2})$, as claimed.

\medskip

Finally we prove Proposition 1.1. Consider $F_n \leq m \leq F_n +
F_{n-3} - 2$. The parts of a partition of $m$ come from the list
$F_2, F_3, \ldots, F_n$. To each partition $\pi$ of $m$, assign
the partition $\pi'$ of $m' = F_{n+2} - 2 - m$ consisting of all
the numbers on the above list that do not appear in $\pi$. Any
partition of $m'$ can be obtained in such a way from a partition
of $m$: the partitions of $m'$ also have all their parts less than
or equal to $F_n$, because it is easily seen that $m' < F_{n+1}$.

So the partitions of $m$ are in bijection with the partitions of
$m'$. If a partition $\pi$ of $m$ has $k$ parts, the
corresponding partition $\pi'$ of $m'$ has $n-1-k$ parts.
Therefore, if $n$ is odd, the bijection takes odd partitions to
odd partitions and even partitions to even partitions, and
$a(m)=a(m')$. If $n$ is even, the bijection takes odd partitions
to even partitions, and even partitions to odd partitions, and
$a(m)=-a(m')$. In any case, $a(m) = (-1)^{n-1} a(m')$.

Now, it is easily seen that $F_n+F_{n-2} \leq m' \leq F_{n+1}-2$.
Therefore Proposition 1.3 applies, and
$a(m')=a(m'-F_n-F_{n-2})=a(F_n+F_{n-3}-2-m)$. Hence $a(m) =
(-1)^{n-1}a(F_n+F_{n-3}-2-m)$, which is what we wanted to show.

\medskip

{\bf Proposition 2.} Given an integer $n$, pick an integer $m$
uniformly at random from the interval $[0,n]$. Let $p_n$ be the
probability that $a(m)=0$ or, equivalently, that $r_E(m)=r_O(m)$.

Then $\lim_{n \rightarrow \infty} p_n = 1$.

\medskip

{\bf Proof.} Let $\alpha_n$ be the number of non-zero
coefficients among the first $F_n$ coefficients $a(0), a(1),
\ldots, a(F_n-1)$. Then $p_{(F_n-1)} = 1 - \alpha_n / F_n$. We
shall prove that $\lim_{n \rightarrow \infty} \alpha_n / F_n =
0$. Proposition 2 follows easily from this.

First we obtain a recurrence relation for $\alpha_n$. Consider
the non-zero coefficients $a(m)$ for $F_n \leq m \leq F_{n+1}-1$.
We know that there are $\alpha_{n+1}-\alpha_n$ such coefficients.
Now split the interval $[F_n, F_{n+1}-1]$ into the three
subintervals $[F_n, F_n+F_{n-3}-2], [F_n+F_{n-3}-1,
F_n+F_{n-2}-1]$ and $[F_n+F_{n-2}, F_{n+1}-1]$. Proposition 1.2
shows that that there are no non-zero coefficients in the second
subinterval, and Proposition 1.3 shows that there are
$\alpha_{n-3}$ non-zero coefficients in the third subinterval.
Because $a(F_{n-3}-1)$ is non-zero for all $n \geq 5$ (this
follows inductively from Proposition 1.3), Proposition 1.1 shows
that there are $\alpha_{n-3}-1$ non-zero coefficients in the first
subinterval. We conclude that $\alpha_{n+1}-\alpha_n =
2\alpha_{n-3}-1$.

The characteristic polynomial of this recurrence relation is
$x^4-x^3-2=0$, and its roots are approximately $r_1 \approx 1.54,
r_2 = -1, r_3 \approx 0.23 + 1.12i$ and $r_4 \approx 0.23 -
1.12i$. It follows from standard results on linear recurrences
that $\alpha_n = O(r_1^{\, n})$, while $F_n = \Theta(\lambda^n)$,
where $\lambda = (\sqrt5 + 1)/2 \approx 1.62$. Since $r_1 <
\lambda$, we conclude that $\lim_{n \rightarrow \infty } \alpha_n
/ F_n = 0$.

\medskip

{\bf Acknowledgement.} The author would like to thank Richard
Stanley for encouraging him to work on this problem, and for
pointing out \cite{4}.

\medskip


\medskip

\begin{thebibliography}{1}

\bibitem{1}
L. Carlitz. ``Fibonacci Representations." {\it The Fibonacci
Quarterly} {\bf 6.4} (1968): 193-220.

\bibitem{2}
H. H. Ferns. ``On the Representations of Integers as Sums
of Distinct Fibonacci Numbers." {\it The Fibonacci Quarterly} {\bf
3.1} (1965): 21-30.

\bibitem{3}
D. Klarner. ``Partitions of $N$ into Distinct Fibonacci
Numbers." {\it The Fibonacci Quarterly} {\bf 6.4} (1968): 235-243.

\bibitem{4}
N. Robbins. ``Fibonacci Partitions." {\it The Fibonacci Quarterly}
{\bf 34.4} (1996): 306-313.

\end{thebibliography}
\end{document}